# A METHOD FOR EXTRACTION OF ARCS OF THE ALGEBRAIC CURVES

B. Malešević [1]
M. Petrović [2]
B. Banjac [3]
I. Jovović [4]
P. Jovanović [5]

**Abstract**

*In this paper we present a method for extraction of arcs of the algebraic curves of the higher order. Method is applied on conics, Cartesian ovals, trifocal curves and generalized Weber's curve.*

---

[1] PhD Branko Malešević, associate professor at Faculty of Electrical Engineering, Department of Applied Mathematics, University of Belgrade, Serbia, e-mail: malesevic@etf.rs

[2] MSc Maja Petrović, assistant at Faculty of Transport and Traffic Engineering, University of Belgrade, Serbia, e-mail: majapet@sf.bg.ac.rs

[3] MSc Bojan Banjac, student of doctoral studies of Software Engineering, Faculty of Electrical Engineering, University of Belgrade, Serbia,
assistant at Faculty of Technical Sciences, Computer Graphics Chair, University of Novi Sad, Serbia, e-mail: bojan.banjac@uns.ac.rs

[4] PhD Ivana Jovović, assistant professor at Faculty of Electrical Engineering, Department of Applied Mathematics, University of Belgrade, Serbia, e-mail: ivana@etf.rs

[5] MSc Petar Jovanović, student of doctoral studies of Computing, Department od Service and Information System Engineering, Universitat Politècnica de Catalunya, BarcelonaTech, Barcelona, Spain, e-mail: petar@essi.upc.edu

## 1. INTRODUCTION

For given three foci (points) and three directrices (lines) in Euclidean plane, we consider the locus of points that satisfy the following equation:

$$\alpha_1 R_1 + \alpha_2 R_2 + \alpha_3 R_3 + \beta_1 r_1 + \beta_2 r_2 + \beta_3 r_3 = \mathcal{S} \qquad (1)$$

where $R_i$ is Euclidean distance from the point of locus to the $i^{th}$ focus, $r_j$ is Euclidean distance from the point of locus to the $j^{th}$ directrix and $\alpha_i, \beta_j, \mathcal{S} \epsilon \mathbb{R}$ $(i,j = 1..3)$. The main objective of this paper is to extract the conditions that enable correct algebraic transformations, which lead either to arcs or complete algebraic curves. The conditions formed in that manner determine the parts of the plane constrained by the algebraic equations of equal or lower order than the algebraic equations of type (1).

The curves satisfying equation (1) are generalization of Weber's trifocal curves [7], (see also [6]), and trough this paper we will call them *generalized Weber's curves*. Generalized Weber's curves allow also introduction of distance to directrix.

## 2. SOME EXAMPLES OF THE GENERALIZES WEBER'S CURVES

### 2.1. Cone sections

**2.1.1. Parabola.** Let the directrix be the line $x = -\frac{p}{2}$ and let the focus be the point $F\left(\frac{p}{2}, 0\right)$, for the focal parameter $p > 0$. Parabola is the locus of points that have the same distance from the directrix as the focus [8]; in other words:

$$R_1 - r_1 = 0, \qquad (2)$$

where $R_1$ is distance to focus and $r_1$ distance to directrix. Equality (2) is equivalent to:

$$\sqrt{\left(x - \frac{p}{2}\right)^2 + y^2} = \left|x + \frac{p}{2}\right|. \qquad (3)$$

Let us notice that $\sqrt{\left(x - \frac{p}{2}\right)^2 + y^2} \geq \left|x - \frac{p}{2}\right| = -x + \frac{p}{2} > -x - \frac{p}{2} = \left|x + \frac{p}{2}\right|$ for $x < -\frac{p}{2}$. Therefore, the equality (3) can be considered only

for the half-plane $x \geq -\frac{p}{2}$, and in that case, by squaring the equality $\sqrt{\left(x-\frac{p}{2}\right)^2 + y^2} = x + \frac{p}{2}$, we obtain the algebraic equation for parabola

$$y^2 = 2px. \tag{4}$$

Previously determined parabola is completely inside the half-plane $x \geq -\frac{p}{2}$ (*Figure 1*).

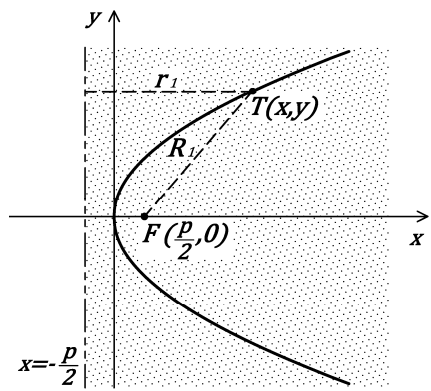

*Figure 1. Parabola*

**2.1.2. Ellipse.** Let $F_1(-c,0)$ and $F_2(c,0)$ be the foci, where $c = \sqrt{a^2 - b^2}$ for $a \geq b > 0$. The ellipse is the locus of points such that the sum of the distances to two foci is constant [8]; in other words:

$$\boldsymbol{R_1 + R_2 = S}, \tag{5}$$

where $R_1$ and $R_2$ are the distances to foci and $S = 2a$. Equality (5) is equivalent to:

$$\sqrt{(x+c)^2 + y^2} + \sqrt{(x-c)^2 + y^2} = 2a. \tag{6}$$

By squaring equality (6) we obtain

$$2\sqrt{((x+c)^2 + y^2)((x-c)^2 + y^2)} = 4a^2 - 2c^2 - 2x^2 - 2y^2. \tag{7}$$

The squaring is correct if $4a^2 - 2c^2 - 2x^2 - 2y^2 \geq 0$, i.e.

$$x^2 + y^2 \leq R^2, \tag{8}$$

where $R = \sqrt{2a^2 - c^2} = \sqrt{a^2 + b^2}$ is a radius of the circumscribed circle of a tangent rectangle (*Figure 2*). All points of ellipse are in the interior area of the circle $x^2 + y^2 = R^2$. By the second sqaring the equality (7) we obtain the canonical form of ellipse:

$$\frac{x^2}{a^2} + \frac{y^2}{b^2} = 1. \tag{9}$$

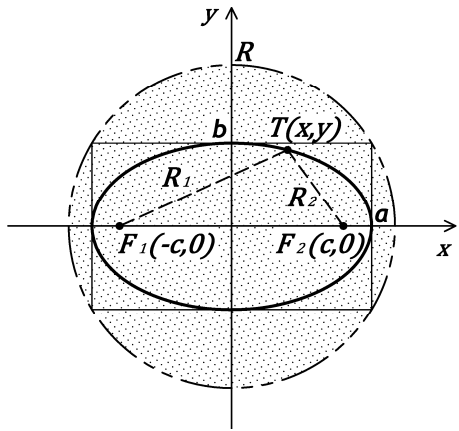

*Figure 2*. Ellipse

**2.1.3. Hyperbola.** Let $F_1(-c, 0)$ and $F_2(c, 0)$ be the foci, where $c = \sqrt{a^2 + b^2}$ for $a \geq b > 0$. The hyperbola is the locus of points such that the absolute value of difference of the distances to two foci is constant [8]; in other words:

$$|R_1 - R_2| = S, \tag{10}$$

where $R_1$ and $R_2$ are the distances to foci and $S = 2a$. Equality (10) is equivalent to:

$$\left|\sqrt{(x+c)^2 + y^2} - \sqrt{(x-c)^2 + y^2}\right| = 2a. \tag{11}$$

By squaring equality (11) we obtain

$$2\sqrt{((x+c)^2 + y^2)((x-c)^2 + y^2)} = -4a^2 + 2c^2 + 2x^2 + 2y^2. \tag{12}$$

The squaring is correct if $-4a^2 + 2c^2 + 2x^2 + 2y^2 \geq 0$, i.e.

$$x^2 + y^2 \geq R^2, \tag{13}$$

where $R = \sqrt{2a^2 - c^2} = \sqrt{a^2 - b^2}$. Radius $R$ of the circle $x^2 + y^2 = R^2$ can be constructively determined from the tangent rectangle and $R < a$ holds (*Figure 3*). All points of hyperbola are in the exterior area of the circle $x^2 + y^2 = R^2$. By the second sqaring the equality (12) we obtain the canonical form of hyperbola:

$$\frac{x^2}{a^2} - \frac{y^2}{b^2} = 1. \tag{14}$$

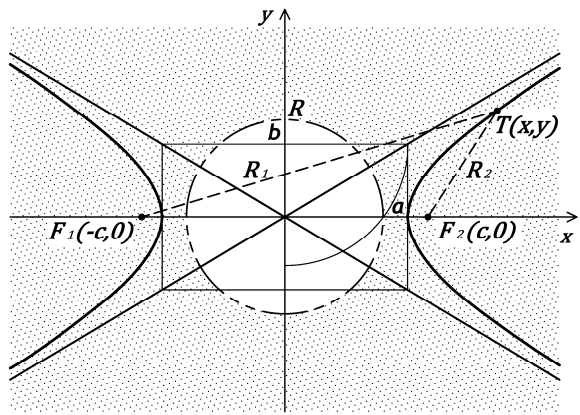

*Figure 3*. Hyperbola

### 2.2. Cartesian Ovals

Let $F_1(-c, 0)$ and $F_2(c, 0)$ be the foci, for $c > 0$. Cartesian ovals is locus of point satisfying:

$$mR_1 \pm nR_2 = \mathcal{S}, \tag{15}$$

where $R_1$ and $R_2$ are distances to foci, $\mathcal{S} > 0$ and $m, n \in \mathbb{N}$, [8]. Equality (15) states:

$$m\sqrt{(x-c)^2 + y^2} \pm n\sqrt{(x+c)^2 + y^2} = \mathcal{S}. \tag{16}$$

By squaring previous equality we obtain

$$\pm 2mn\sqrt{((x-c)^2 + y^2)((x+c)^2 + y^2)} = \mathcal{S}^2 - (m^2 + n^2)(x^2 + y^2)$$
$$+ 2c(m^2 - n^2)x - c^2(m^2 + n^2). \tag{17}$$

Let us form a circle

$$\left(x - \frac{c(m^2-n^2)}{m^2+n^2}\right)^2 + y^2 = R^2 \qquad (18)$$

with radius $R = \frac{1}{m^2+n^2}\sqrt{S^2(m^2+n^2) - 4c^2m^2n^2} > 0$ for parameter $S > c\frac{2mn}{\sqrt{m^2+n^2}}$. There is two possible cases.

1. Let us consider the operand plus (+) in formula (15), i.e. a Cartesian oval $mR_1 + nR_2 = S$. Then, the left side of equality (17) is positive and further squaring is correct if $S^2 - (m^2+n^2)(x^2+y^2) + 2c(m^2-n^2)x - c^2(m^2+n^2) \geq 0$, i.e. for the interior points of the circle

$$\left(x - \frac{c(m^2-n^2)}{m^2+n^2}\right)^2 + y^2 \leq R^2. \qquad (19)$$

2. Let us consider the operand plus (−) in formula (15), i.e. a Cartesian oval $mR_1 - nR_2 = S$. Then, the left side of equality (17) is negative and further squaring is correct if $S^2 - (m^2+n^2)(x^2+y^2) + 2c(m^2-n^2)x - c^2(m^2+n^2) \leq 0$, i.e. for the exterior points of the circle

$$\left(x - \frac{c(m^2-n^2)}{m^2+n^2}\right)^2 + y^2 \geq R^2. \qquad (20)$$

In the first case, the Cartesian oval is in the interior area of the circle, while in the second case the Cartesian oval is in the exterior area of the circle (*Figure 4*). In both cases, by squaring equality (17) we obtain the following algebraic equation:

$$(m^2-n^2)^2 x^4 + 2(m^2-n^2)x^2y^2 + (m^2-n^2)^2 y^4 \qquad (21)$$
$$-4c(m^4-n^4)x^3 - 4c(m^4-n^4)xy^2$$
$$+(2c^2(3m^4 + 2m^2n^2 + 3n^2) + 2S^2(m^2+n^2))x^2$$
$$+(2c^2(n^4 - 2m^2n^2 + n^4) + 2S^2(m^2+n^2))y^2$$
$$+4c(m^2-n^2)(S^2 - c^2(m^2+n^2))x$$
$$+(S + c(m+n))(S + c(m-n))(S + c(-m+n))(S + c(-m-n)) = 0.$$

As a special case, if $m = n$, then the Cartesian oval transforms into ellipse; while if $m = -n$, then the Cartesian oval transforms into hyperbola.

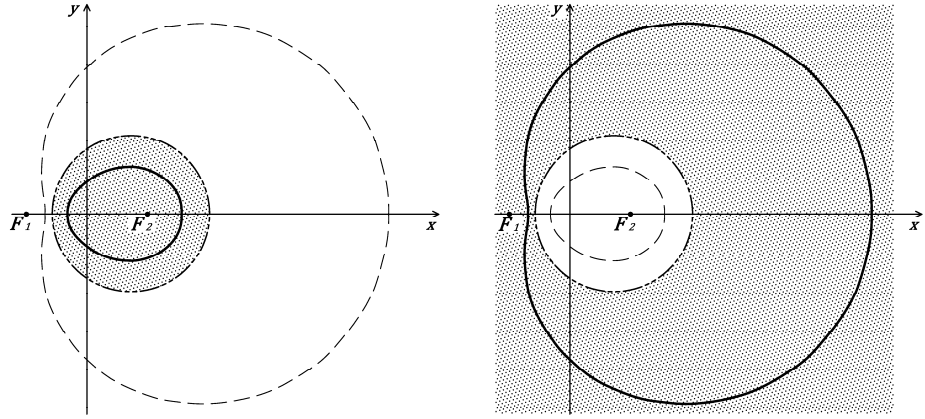

*Figure 4. Cartesian ovals*

R. Descartes had obtained these ovals during his research related to optics, i.e. they resulted as a solution to the problem of searching for a curve such that the rays stemming from one point refract on this curve and afterwards pass through another given point. Light rays stemming from the focus $F_1$, after refracting in the points of the Cartesian oval, meet in the focus $F_2$. In other words, these ovals share two mediums with the refractive index $n/m$. M. Chasles had further studied these ovals, and showed that Cartesian ovals also have the third focus, i.e. that the following holds:

$$mR_1 + nR_2 = S, \quad m_1R_1 + n_1R_3 = S_1, \quad m_2R_2 + n_2R_3 = S_2.$$

By multiplying the first equality by $S_1$ and the second one by $S$ and by subtracting he had obtained:

$$\alpha R_1 + \beta R_2 + \gamma R_3 = 0 \tag{22}$$

where $\alpha = mS_1 - m_1S$, $\beta = nS_1$, $\gamma = n_1S \in \mathbb{R}$. Based on the previous equality, one can conclude that Cartesian ovals are loci of points for which the sum of distances from three fixed points, previously multiplied by some numerical values, is annulled.

### 2.3. Trifocal ellipse

The standard definition of a trifocal ellipse is given by the following equation:

$$R_1 + R_2 + R_3 = S \tag{23}$$

where $R_1 = \sqrt{(x-p)^2 + y^2}$, $R_2 = \sqrt{(x-q)^2 + y^2}$, and $R_3 = \sqrt{x^2 + (y-r)^2}$ are Euclidean distances of the point $(x, y)$ to three foci $F_1(p, 0)$, $F_2(q, 0)$, $F_3(0, r)$, respectively, and for a given parameter $S > 0$, $(p, q, r \in \mathbb{R})$ [2], [6], [7]. The Fermat-Torricelli point determines the minimal value of parameter $S = S_0$. If $0 < S < S_0$ then a trifocal ellipse does not exist, and for $S > S_0$ trifocal the ellipse is a non-degenerated egg curve (*Figure 5*).

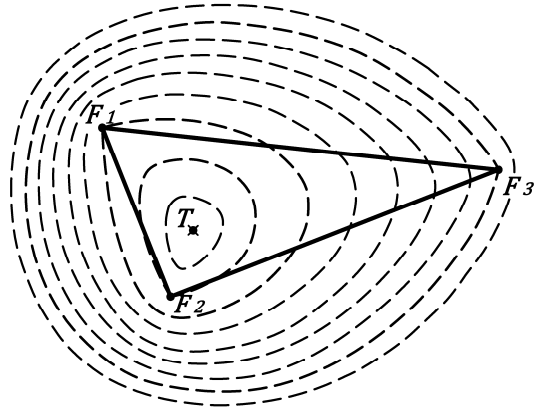

**Figure 5.** Trifocal curves

The subject of this section are the trifocal ellipses (23). Let there be given expressions: $Q_1 = (x-p)^2 + y^2$, $Q_2 = (x-q)^2 + y^2$, $Q_3 = x^2 + (y-r)^2$ and $S = const$, $S > S_0$. Let us consider the following algebraic transformation of equation (23) $\Leftrightarrow \sqrt{Q_1} + \sqrt{Q_2} + \sqrt{Q_3} = S$:

$$\sqrt{Q_1} + \sqrt{Q_2} = S - \sqrt{Q_3}. \tag{24}$$

The first squaring requires that the condition $S \geq \sqrt{Q_3}$ holds, which can be reduced to

$$x^2 + (y-r)^2 \leq \rho_1^2 \wedge \rho_1 = S. \tag{25}$$

For the interior points of the circle (25), by squaring the equality (24) we obtain

$$2\sqrt{Q_1}\sqrt{Q_2} + 2S\sqrt{Q_3} = S^2 + Q_3 - Q_2 - Q_1. \tag{26}$$

The second squaring requires that the condition $S^2 \geq Q_1 + Q_2 - Q_3$ holds, which can be reduced to

$$(x - p - q)^2 + (y + r)^2 \leq \rho_2^2 \wedge \rho_2 = \sqrt{S^2 + 2pq + 2r^2}, \quad (27)$$

for $S > 0$ such that $S^2 + 2pq + 2r^2 > 0$. Then, for the points within the intersection of interiors of the circles (25) and (27), by squaring the equality (26) we obtain

$$8S\sqrt{Q_1}\sqrt{Q_2}\sqrt{Q_3} = S^4 - 2S^2Q_3 - 2S^2Q_2 - 2S^2Q_1 \quad (28)$$
$$-2Q_3Q_2 - 2Q_3Q_1 - 2Q_1Q_2 + Q_1^2 + Q_2^2 + Q_3^2$$

The third squaring requires that the condition $S^4 - 2S^2Q_3 - 2S^2Q_2 - 2S^2Q_1 - 2Q_3Q_2 - 2Q_3Q_1 - 2Q_1Q_2 + Q_1^2 + Q_2^2 + Q_3^2 \geq 0$ holds, which can be reduced to

$$\alpha = 3x^4 + 6x^2y^2 + 3y^4 - 4(p+q)x^3 - 4rx^2y - 4(p+q)xy^2 - 4ry^3 \quad (29)$$
$$+2(3S^2 + r^2 - p^2 + 4pq - q^2)x^2 + 8r(p+q)xy + 2(3S^2 + p^2 + q^2 - r^2)y^2$$
$$-4(p+q)(S^2 + r^2 - (p-q)^2)x - 4r(S^2 + p^2 + q^2 - r^2)y$$
$$-S^4 + 2(p^2 + q^2 + r^2)S^2 + (p+q+r)(p+q-r)(p-q+r)(-p+q+r) \leq 0$$

The previous fourth order curve $\alpha = 0$ can be written in the following form

$$\alpha = k(x^2 + y^2)^2 + (Ax + By)(x^2 + y^2) \quad (30)$$
$$+Cx^2 + Dxy + Ey^2 + Fx + Gy + H = 0$$

and it belongs to the family of the curves from [3] (see formula (3.5)). Considering the points within the intersection of interiors of the circles (25) and (27) with the interior of the curve (29), by squaring the equation (28) we finally obtain the eighth order curve:

$$64S^2Q_1Q_2Q_3 - ((S^2 + Q_3 - Q_2 - Q_1)^2 - 4Q_1Q_2 - 4S^2Q_3)^2 = 0 \quad (31)$$

which can be explicitly written as

$$64S^2((x-p)^2 + y^2)((x-q)^2 + y^2)(x^2 + (y-r)^2) \quad (32)$$
$$-((S^2 + x^2 + (y-r)^2 - (x-q)^2 - (x-p)^2 - 2y^2)^2$$
$$- 4((x-p)^2 + y^2)((x-q)^2 + y^2) - 4S^2(x^2 + (y-r)^2))^2 = 0$$

Let us notice that in the previous derivation we have obtained a sequence of conditions that determine the adjacent areas of the two interiors of the corresponding circles and the interior of the fourth order curve. In the intersection of these three interior areas, we extract, step by step, the part of the algebraic curve (32) that contains closed contour of a trifocal ellipse (*Figure 6*).

Algebraic equation for a trifocal ellipse is derived by means of the elementary algebraic transformations. The obtained equation form coincides with the equation form obtained by using determinants [5] and by means of applying Gröbner basis [2]. It should be noted that the result in the expanded determinant form of the trifocal ellipse representation consists of 2355 sumands [5].

In this section, we primarily specify the exact part of the algebraic curve that represents a trifocal ellipse, and the part of the curve that represents so-called Zarinski closure of the trifocal curve, within the context of the algebraic geometry [5].

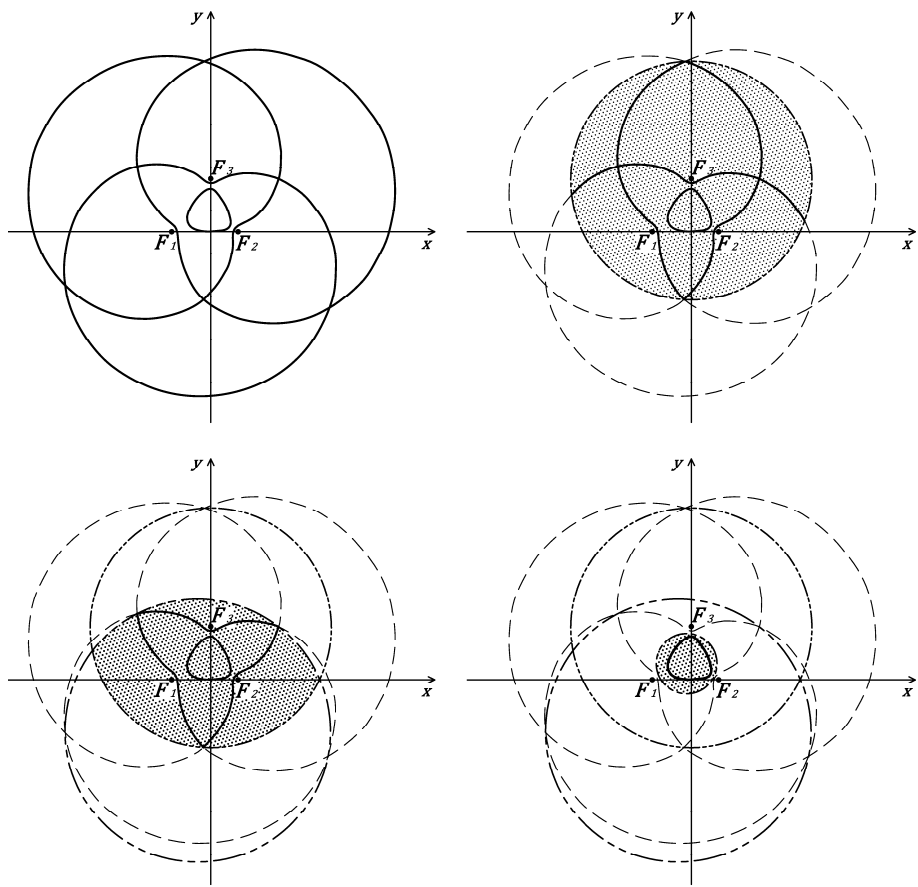

*Figure 6. Algebraic curve and areas for extracting the curves*

## 3. THE ARCHES OF GENERALIZED WEBER'S CURVES

In this section, we consider a generalized Weber's curve (1). One example of such curve is Erdös-Mordell's curve defined by equality:

$$w(x,y) = R_1 + R_2 + R_3 - 2r_1 - 2r_2 - 2r_3 = 0 \qquad (33)$$

which has been introduced and studied in [1], [4]. Considered curve (*Figure 7*) consists of six arches of algebraic curves of the eight order which can be extracted in similar manner as it was given in the previous sections.

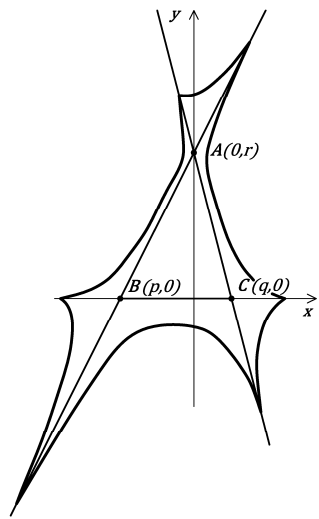

*Figure 7. Erdös – Mordell's curve*

The interior of Erdös-Mordell's curve is determined with $w(x,y) \leq 0$ and inside of it, the Erdös-Mordell's inequality holds:

$$R_1 + R_2 + R_3 \geq 2r_1 + 2r_2 + 2r_3. \qquad (34)$$

P. Erdös in 1935. has set the hypothesis that the inequality (34) holds in the interior of the triangle ABC (where A, B, C are given foci). The hypothesis has been proven in 1937. by L.J. Mordell and D.F. Barrow. In the paper [2] the Erdös-Mordell's inequality has been extended from the interior of the triangle to the interior of the Erdös-Mordell's curve.

In equality (1), the absolute values appear in the distances $r_1, r_2, r_3$. By freeing from the absolute values in some of the regions of a plane, the equality (1) is transformed into the following equalities

$$\alpha_1 R_1 + \alpha_2 R_2 + \alpha_3 R_3 = A_i x + B_i y + C_i \qquad (35)$$

for corresponding $A_i, B_i, C_i \in \mathbb{R}$ $(i = 1..8)$. Using the procedure as in the Section 2.3, some of the arches of the generalized Weber's curve can be determined in the intersection of the corresponding areas, constrained by the appropriate algebraic curves of the second and fourth order. Generalized Weber's curve consists of arches of algebraic curves and these algebraic curves are obtained from (35) with most three squaring.

*Acknowledgement*. Research is partially supported by the Ministry of Science and Education of the Republic of Serbia, Grant No. III 44006 and ON 174032.